
\baselineskip=14pt
\parskip=10pt
\def\Tilde{\char126\relax}
\def\halmos{\hbox{\vrule height0.15cm width0.01cm\vbox{\hrule height
 0.01cm width0.2cm \vskip0.15cm \hrule height 0.01cm width0.2cm}\vrule
 height0.15cm width 0.01cm}}
\font\eightrm=cmr8  
\font\eighttt=cmtt8
\magnification=\magstephalf

\parindent=0pt
\overfullrule=0in
\bf
\centerline
{A CONDENSED CONDENSATION PROOF OF A DETERMINANT EVALUATION}
\centerline
{CONJECTURED BY Greg KUPERBERG AND Jim PROPP}
\rm
\bigskip
\centerline{ {\it Tewodros AMDEBERHAN}{$\,^1$} and 
{\it Shalosh B. EKHAD}\footnote{$^1$}
{\eightrm  \raggedright
Department of Mathematics, Temple University,
Philadelphia, PA 19122, USA. 
{\eighttt E-mail:[tewodros,ekhad]@math.temple.edu; 
WWW: http://www.math.temple.edu/\Tilde [tewodros,ekhad] .}
} 
}

Greg Kuperberg and Jim Propp [P] have conjectured the following determinant 
identity:
$$
\det \left [ {{i+j} \choose {i}}
{{2n-i-j} \choose {n-i}}_{0 \leq i,j \leq n}\right ]
=
 {{ (2n+1)!^{n+1}} \over
 {(2n+1)!!}
} \quad ,
$$
where $a!!:=0! \cdot 1!\cdot 2! \cdots a!$, and $a!:=1\cdot 2 \cdots a$.  
 
This is the special case ($m=n, a=b=0$) of 
$$
\det \left [ {{i+j+a+b} \choose {i+a}}{{2n-i-j-a-b} \choose {n-i-a}}
_{0 \leq i,j \leq m}\right ]
=
$$
$$
{{ (a+b)! (2n+1)!^{m+1} (2n-m)!!m!!(m+a+b)!!(2n-m-a-b)!!a!!b!!
(n-m-a-1)!!(n-m-b-1)!!} \over
 {a!b! (2n+1)!! (n-a)!!(n-b)!!(m+a)!!(m+b)!!(a+b)!!(2n-2m-a-b-1)!!}
} \quad ,
\eqno(Rabbit)
$$
which follows immediately from Dodgson's[D] rule for
evaluating determinants: (For any $n \times n$ matrix $A$,
let $A_r(k,l)$ be the $r \times r$ connected submatrix whose upper leftmost
corner is the entry $a_{k,l}$,)
$$
\det A \,\,=\,\,
{{
\det A_{n-1}(1,1) \det A_{n-1}(2,2) 
-
\det A_{n-1}(1,2) \det A_{n-1}(2,1) }
\over
{\det A_{n-2}(2,2)}} \quad .
\eqno(Lewis)
$$
Indeed, let the left and right sides of $(Rabbit)$
be $L_m (a,b)$ and $R_m (a,b)$
respectively. Dodgson's rule immediately implies that the recurrence:
$$
X_m(a,b)={{X_{m-1}(a,b) X_{m-1}(a+1,b+1)-X_{m-1} (a+1,b) X_{m-1} (a,b+1)}
\over {X_{m-2} (a+1,b+1)}} \quad,
$$
holds with $X=L$.
Since $L_m (a,b)=R_m (a,b)$ for $m=0,1$ (check!), and the recurrence
also holds with $X=R$ (check!\footnote{$^2$}
{\eightrm Divide both sides by the left, then use r!!/(r-1)!!=r!
whenever possible, and then r!/(r-1)!=r whenever possible, reducing
it to a completely routine polynomial identity. The small Maple
package {\eighttt rabbit} obtainable from our Home Pages, performs these
steps mechanically.}
), 
it follows by induction that $L_m(a,b)=R_m(a,b)$ for {\it all} $m$. \halmos.
 
The present proof is in the spirit of [Z1].
Another proof can be found in [A1]. 
The same method yields a q-analog of $(Rabbit)$, that can be found in [A2]. 
A beautiful combinatorial proof of $(Lewis)$ can be found in [Z2].
An alternative proof of $(Rabbit)$ is given in [K].
\eject
{\bf References}
 
[A1] T. Amdeberhan, {\it A WZ proof of a determinant evaluation conjectured by
Kuperberg and Propp}, exclusively published in Amdeberhan's
Home Page {\tt http://www.math.temple.edu/\Tilde tewodros}.
 
[A2] T. Amdeberhan, {\it A q-generalization of a determinant evaluation
conjectured by
Kuperberg and Propp}, exclusively published in Amdeberhan's
Home Page {\tt http://www.math.temple.edu/\Tilde tewodros}.
 
[D] C.L. Dodgson, {\it Condensation of Determinants}, Proceedings
of the Royal Society of London {\bf 15}(1866), 150-155.

[K] C. Krattenthaler, {\it E-mail message to T. Amdeberhan}, dated
9 Aug. 1996, $19:22:27 +0100$ (MET). 
By kind permission of Krattenthaler, it can be gleaned
at the first author's Home Page 
{\tt http://www.math.temple.edu/\Tilde tewodros}.
 
[P] J. Propp, {\it E-mail message to D. Zeilberger}, dated
1 July 1996, 15:47:56 (EDT). 
By kind permission of Propp, it can be gleaned
at the first author's Home Page 
{\tt http://www.math.temple.edu/\Tilde tewodros}.
  
[Z1] D. Zeilberger, {\it Reverend Charles to the aid of Major Percy and
Fields Medalist Enrico}, Amer. Math. Monthly {\bf 103}(1996),
501-502.
 
[Z2] D. Zeilberger, 
{\it Dodgson's determinant-evaluation rule proved by 
Two-Timing Men and Women}, to appear in the Elec. J. of
Combinatorics [Wilf Festschrifft volume]. It can be downloaded
from Zeilberger's Home Page:
{\tt http://www.math.temple.edu/\Tilde zeilberg}.

Version of Sept. 4, 1996. First version: Aug. 2, 1996. 
\bye